\documentclass{amsart}

\usepackage{amssymb,bm,mathrsfs,enumerate,color,scalerel}
\usepackage[noadjust]{cite}
\usepackage{tikz}
\usetikzlibrary{positioning,calc,graphs,quotes,graphs.standard,arrows.meta}
\usepackage{hyperref}

\usetikzlibrary{svg.path}
\definecolor{orcid_color}{HTML}{A6CE39}

\hypersetup{pdfborder={0 0 0}} 
\DeclareRobustCommand{\orcidicon}{%
	\raisebox{.2mm}{\scalerel*{%
	\begin{tikzpicture}[xscale=1,yscale=-1,transform shape]
	\filldraw[color=orcid_color] svg {M256,128c0,70.7-57.3,128-128,128C57.3,256,0,198.7,0,128C0,57.3,57.3,0,128,0C198.7,0,256,57.3,256,128z};
	\filldraw[color=white] svg {M86.3,186.2H70.9V79.1h15.4v48.4V186.2z} svg {M108.9,79.1h41.6c39.6,0,57,28.3,57,53.6c0,27.5-21.5,53.6-56.8,53.6h-41.8V79.1z M124.3,172.4h24.5
		c34.9,0,42.9-26.5,42.9-39.7c0-21.5-13.7-39.7-43.7-39.7h-23.7V172.4z} svg {M88.7,56.8c0,5.5-4.5,10.1-10.1,10.1c-5.6,0-10.1-4.6-10.1-10.1c0-5.6,4.5-10.1,10.1-10.1
		C84.2,46.7,88.7,51.3,88.7,56.8z};
	\end{tikzpicture}}{|}}%
}
\newcommand{\orcid}[1]{\href{https://orcid.org/#1}{\orcidicon}}

\theoremstyle{plain}
\newtheorem{thm}{Theorem}[section]
\theoremstyle{definition}
\newtheorem{exam}[thm]{Example}
\newtheorem{conj}[thm]{Conjecture}

\newcommand{\arxiv}[1]{arXiv:\href{http://arxiv.org/abs/#1}{#1}}

\begin{document}

\title{A duality of scaffolds for translation association schemes}
\author{Xiaoye Liang}
\address{School of Mathematical Sciences, Anhui University, Hefei, Anhui 230601, PR China}
\email{liangxy0105@foxmail.com}
\author[Ying-Ying Tan]{Ying-Ying Tan\,\orcid{0000-0001-8152-6722}}
\address{School of Mathematics and Physics, Anhui Jianzhu University, Hefei, Anhui 230601, PR China}
\email{tansusan1@ahjzu.edu.cn}
\thanks{Y.-Y Tan is supported by the National Natural Science Foundation of China (No. 11801007 and No. 12171002).}
\author[Hajime Tanaka]{Hajime Tanaka\,\orcid{0000-0002-5958-0375}}
\address{\href{http://www.math.is.tohoku.ac.jp/}{Research Center for Pure and Applied Mathematics}, Graduate School of Information Sciences, Tohoku University, Sendai 980-8579, Japan}
\email{htanaka@tohoku.ac.jp}
\thanks{H. Tanaka was supported by JSPS KAKENHI Grant Number JP20K03551.}
\author{Tao Wang}
\address{Qingdao Pennon Education Technology Co., Ltd, Qingdao, Shandong 266101, PR China}
\email{wangtao998877@foxmail.com}

\begin{abstract}
\emph{Scaffolds} are certain tensors arising in the study of association schemes, and have been (implicitly) understood diagrammatically as digraphs with distinguished ``root'' nodes and with matrix edge weights, often taken from Bose--Mesner algebras.
In this paper, we first present a slight modification of Martin's conjecture (2021) concerning a duality of scaffolds whose digraphs are embedded in a closed disk in the plane with root nodes all lying on the boundary circle, and then show that this modified conjecture holds true if we restrict ourselves to the class of \emph{translation} association schemes, i.e., those association schemes that admit abelian regular automorphism groups.
\end{abstract}

\maketitle

\hypersetup{pdfborder={0 0 1}} 

\section{Introduction}

Certain tensors arise in the study of (commutative) association schemes.
One of the first instances appears in the work of Cameron, Goethals, and Seidel \cite{CGS1978IM} on the implications of the (non)vanishing of the Krein parameters of association schemes.
Initiated by Neumaier around 1989, several researchers have used these tensors to prove results about association schemes and related objects; see, e.g., \cite{Dickie1995D,Jaeger1995JAC,PN2021+a,PN2021+b,Suzuki1998JACa,Suzuki1998JACb,Terwilliger1987JCTA}.
These tensors have been (implicitly) understood diagrammatically.
For example, those appearing in \cite{CGS1978IM} are $3$-tensors and are  encoded as
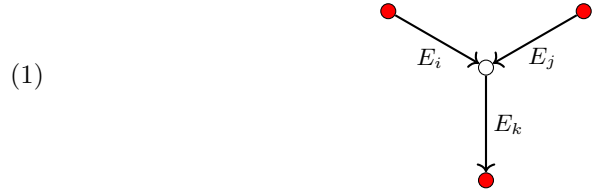
\begin{equation}\label{star scaffold}
\begin{array}{c}
\begin{tikzpicture}
\draw[magenta!0] (0,0) circle [radius=1.5];
\graph [ empty nodes, nodes={circle,draw,fill=red,inner sep=2pt}, edges={thick}]  { subgraph I_n [n=3, clockwise,radius=15mm,phase=30,V={a,b,c}];  {a,c} -> d[fill=white] -> b };
\node at (.75,.12) {\small $E_j$};
\node at (.3,-.75) {\small $E_k$};
\node at (-.75,.12) {\small $E_i$};
\end{tikzpicture}
\end{array}
\end{equation}
where $E_i,E_j$, and $E_k$ are primitive idempotents of an association scheme.
Martin \cite{Martin2021LAA} called these tensors (and also their diagrams) \emph{scaffolds}, and studied them systematically and comprehensively.
The red nodes are referred to as the \emph{root nodes} or simply \emph{roots}.
In general, the order of such a tensor equals the number of roots in its diagram.
We allow the case where there is no root; scaffolds are then $0$-tensors, i.e., scalars, and are precisely the \emph{partition functions} as defined by Jaeger \cite{Jaeger1995JAC}.

For most of the scaffolds discussed in the literature, the underlying digraphs may be given embeddings in a closed disk in the plane, with root nodes all lying on the boundary circle, and the edge weights consist only of the adjacency matrices $A_i$ or of the primitive idempotents $E_i$.
In view of the dual roles played by the $A_i$ and the $E_i$, we then have an apparent duality for such a scaffold; that is to say, we take the faces in the interior of the disk as the nodes of the dual, where two faces are joined when they are separated by an edge (with direction taken $90^{\circ}$ clockwise), and where the $A_i$ and the $E_i$ are interchanged in the edge weights.\footnote{This, however, is not a duality in the strict sense, in that if we take the dual twice then we return to the original scaffold except that all the edges are reversed.}
For example, the above scaffold is the dual of
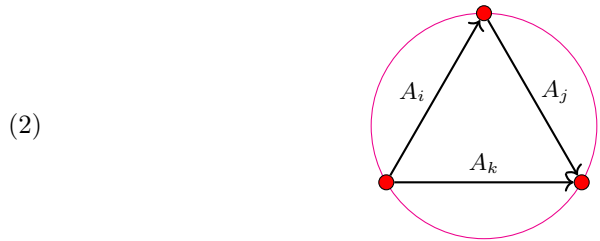
\begin{equation}\label{triangle scaffold}
\begin{array}{c}
\begin{tikzpicture}
\draw[magenta] (0,0) circle [radius=1.5];
\graph [ empty nodes, nodes={circle,draw,fill=red,inner sep=2pt}, edges={thick}]  { subgraph I_n [n=3, clockwise,radius=15mm,V={a,b,c}]; a -> b <- c -> a };
\node at (.95,.45) {\small $A_j$};
\node at (0,-.5) {\small $A_k$};
\node at (-.95,.45) {\small $A_i$};
\end{tikzpicture}
\end{array}
\end{equation}
where we also drew the boundary circle for convenience.
Many combinatorial statements about association schemes were first interpreted as implications among equations of scaffolds whose edge weights are taken from the $A_i$, and then dualized.
Moreover, the validity of these dual implications was usually verified by also dualizing the proofs.
Martin \cite[Conjecture 4.1]{Martin2021LAA} conjectured that this process is always possible.
Namely, his conjecture asserts that if an implication among scaffold equations with edge weights  from the $A_i$ is true for all association schemes, then so is the dual implication.

In this paper, after recalling necessary definitions, we first present in Section \ref{sec: definitions and results} a slightly modified version of Martin's conjecture.
In particular, we make the concept of dual scaffold equations more precise.
See Conjecture \ref{Martin conjecture} and the subsequent comments.
We then show in Theorem \ref{main theorem} that Conjecture \ref{Martin conjecture} holds true if we restrict ourselves to the class of \emph{translation} association schemes \cite[Section 2.10]{BCN1989B}, i.e., those association schemes that admit abelian regular automorphism groups.
For such an association scheme, we can naturally define its dual, which is again a translation association scheme.
Examples in this class include the Hamming schemes, the Doob schemes, and the four forms schemes (bilinear, alternating, Hermitian, and quadratic); cf.~\cite[Section III.6]{BI1984B}.
We obtain Theorem \ref{main theorem} as an immediate consequence of Theorem \ref{real main theorem}, which will be proved in Section \ref{sec: proof}.

\section{Scaffolds, association schemes, and Bose--Mesner algebras}
\label{sec: definitions and results}

For a nonempty finite set  $X$, let $\mathsf{Mat}_X(\mathbb{C})$ be the $\mathbb{C}$-algebra of complex matrices with rows and columns indexed by $X$.
We also let $\mathbb{C}^X$ be the $\mathbb{C}$-vector space of complex column vectors with coordinates indexed by $X$, and let $\{\hat{x}:x\in X\}$ denote the standard basis of $\mathbb{C}^X$.
We always endow $\mathbb{C}^X$ with the standard Hermitian inner product $\langle\hat{x},\hat{y}\rangle=\delta_{x,y}$ $(x,y\in X)$.

For two nonempty finite sets $X$ and $Y$, we let $Y^X$ denote the set of functions from $X$ to $Y$.
For $\gamma\in Y^X$, the image of $x\in X$ will be denoted by $x^{\gamma}$, except when both $X$ and $Y$ are finite groups and $\gamma$ is a homomorphism, in which case we will write $\gamma(x).$

Suppose that we are given
\begin{itemize}
\setlength{\itemsep}{1mm}
\item a finite directed graph $G$ with vertex set $V$ and edge set $E$, possibly with loops and/or multiple edges;
\item an ordered multiset $R=\{r_1,\dots,r_{\ell}\}\subset V$ of \emph{root nodes} or simply \emph{roots};
\item a nonempty finite set $X$ and an edge weight function $w:E\rightarrow \mathsf{Mat}_X(\mathbb{C})$.
\end{itemize}
We denote the head and tail of $e\in E$ by $h(e)$ and $t(e)$, respectively.
Later we will assume that the set $X$ is the vertex set of an association scheme, so we call the elements of $V$ \emph{nodes} to distinguish $V$ from $X$.
The (\emph{symmetric}) \emph{scaffold} $\mathsf{S}(G,R;w)$ \emph{of order} $\ell$ \cite{Martin2021LAA} is the vector in $\mathbb{C}^{X^{\ell}}=(\mathbb{C}^X)^{\otimes\ell}$ defined by
\begin{equation*}
	\mathsf{S}(G,R;w)= \sum_{\sigma\in X^V}\!\left(\prod_{e\in E}w(e)_{t(e)^{\sigma}\!,h(e)^{\sigma}} \!\right) \widehat{\bm{r}^{\sigma}},
\end{equation*}
where we write $\widehat{\bm{r}^{\sigma}}:=\widehat{r_1^{\sigma}} \otimes\cdots\otimes\widehat{r_{\ell}^{\sigma}}$ for brevity.
Martin \cite{Martin2021LAA} also considered the case where $\sigma$ ranges over a certain subset of $X^V$, but we will only discuss the above unrestricted case in this paper.

We represent the scaffold $\mathsf{S}(G,R;w)$ pictorially with a plane drawing (possibly with crossings) of $G$ with root nodes written in red and edges labelled with their weights.
For example, if $G$ is the directed path %
\begin{tikzpicture}
\graph [ empty nodes, nodes={circle,draw,fill=white,inner sep=2pt}, edges={thick}] { a -> b -> c };
\end{tikzpicture},
$R=\{r_1,r_2\}$ where $r_1$ (resp.~$r_2$) is the left (resp.~right) end node, and the weight of the left (resp.~right) edge is $A$ (resp.~$B$), then $\mathsf{S}(G,R;w)$ is represented as
\begin{center}
\begin{tikzpicture}
\graph [ empty nodes, nodes={circle,draw,inner sep=2pt}, edges={thick}] { a[fill=red,label=left:$r_1$] ->["\small $A$"] b[fill=white] ->["\small $B$"] c[fill=red,label=right:$r_2$] };
\end{tikzpicture}.
\end{center}
Working with the diagrams helps manipulating scaffolds, because this allows us to understand and handle various scaffold equations pictorially.
For example, we have
\begin{center}
\begin{tikzpicture}
\graph [ empty nodes, nodes={circle,draw,inner sep=2pt}, edges={thick}] { a[fill=red,label=left:$r_1$] ->["\small $A$"] b[fill=white] ->["\small $B$"] c[fill=red,label=right:$r_2$] };
\node at (3.5,0) {$=$};
\graph [ empty nodes, nodes={xshift=5cm,circle,draw,inner sep=2pt,fill=red}, edges={thick}] { a[label=left:$r_1$] ->["\small $AB$"] b[label=right:$r_2$] };
\end{tikzpicture}.
\end{center}
Another example is
\begin{center}
\begin{tikzpicture}
\graph [ empty nodes, nodes={circle,draw,inner sep=2pt,fill=red}, edges={thick}]
{ a[label=left:$r_1$] ->[bend left=50,"\small $A$"] b[label=right:$r_2$], a ->[bend right=50,"\small $B$"'] b };
\node at (2.5,0) {$=$};
\graph [ empty nodes, nodes={xshift=4cm,circle,draw,inner sep=2pt,fill=red}, edges={thick}]
{ a[label=left:$r_1$] ->["\small $A\circ B$"] b[label=right:$r_2$] };
\end{tikzpicture}
\end{center}
where $\circ$ denotes entrywise (or \emph{Hadamard}) multiplication.
See \cite[Appendix A]{Martin2021LAA} for a list of rules for the manipulation of scaffolds.
We often \emph{identify} these diagrams with the corresponding scaffolds.

We now recall association schemes and their Bose--Mesner algebras.
Let $\mathcal{R}=\{R_0,\dots,R_d\}$ be a set of binary relations on $X$.
For every $i$, let $A_i$ be the adjacency matrix of the (di)graph $(X,R_i)$.
We say that the pair $(X,\mathcal{R})$ is a (\emph{commutative}) \emph{association scheme} \emph{with} $d$ \emph{classes} if
\begin{enumerate}[({AS}1)]
\setlength{\itemsep}{1mm}
\item\label{AS1} $A_0=I$, the identity matrix;
\item\label{AS2} $\sum_{i=0}^dA_i=J$, the all ones matrix;
\item\label{AS3} $A_i^{\mathsf{T}}\in\{A_0,\dots,A_d\}$ for $0\leqslant i\leqslant d$, where $^{\mathsf{T}}$ denotes transpose;
\item\label{AS4} $A_iA_j=A_jA_i\in\bm{A}:=\operatorname{span}_{\mathbb{C}}\{A_0,\dots,A_d\}$ for $0\leqslant i,j\leqslant d$.
\end{enumerate}
By (AS\ref{AS1}), (AS\ref{AS2}), and (AS\ref{AS4}), the vector space $\bm{A}$ is a $(d+1)$-dimensional commutative $\mathbb{C}$-algebra, and is called the \emph{Bose--Mesner} \emph{algebra} of $(X,\mathcal{R})$.
The algebra $\bm{A}$ is semisimple as it is closed under conjugate-transpose by (AS\ref{AS3}), and hence it has another basis consisting of the primitive idempotents $E_0=|X|^{-1}J,E_1,\dots,E_d$, i.e., $E_iE_j=\delta_{i,j}E_i$, $\sum_{i=0}^dE_i=I$; cf.~\cite[Section II.3]{BI1984B}, \cite[Section 2.2]{BCN1989B}.
By (AS\ref{AS2}), $\bm{A}$ is also closed under $\circ$.
The $A_i$ are the primitive idempotents of $\bm{A}$ with respect to $\circ$, i.e., $A_i\circ A_j=\delta_{i,j}A_i$, $\sum_{i=0}^dA_i=J$.
The \emph{first} and \emph{second eigenmatrices} $P$ and $Q$ are defined respectively by
\begin{equation}\label{eigenmatrices}
	A_i=\sum_{j=0}^dP_{j,i}E_j, \qquad E_i=|X|^{-1}\sum_{j=0}^dQ_{j,i}A_j \qquad (0\leqslant i\leqslant d).
\end{equation}

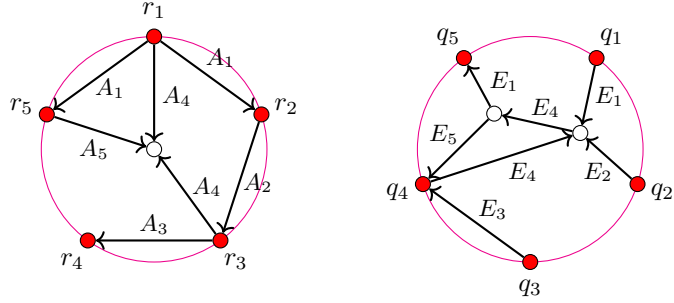
\begin{figure}
\begin{tikzpicture}
\draw[magenta] (0,0) circle [radius=1.5];
\graph [ empty nodes, nodes={circle,draw,fill=red,inner sep=2pt}, edges={thick}]  { subgraph I_n [n=5, clockwise,radius=15mm,V={a,b,c,d,e}];  a -> b -> c -> d, a -> e, {a,c,e} -> f[fill=white] };
\foreach \x in {1,...,5} \node at (90+72-72*\x:1.85) {$r_{\x}$};
\node at (54:1.5) {\small $A_1$};
\node at (342:1.45) {\small $A_2$};
\node at (270:1) {\small $A_3$};
\node at (126:1) {\small $A_1$};
\node at (.3,.7) {\small $A_4$};
\node at (-.8,0) {\small $A_5$};
\node at (.7,-.5) {\small $A_4$};
\draw[magenta] (5,0) circle [radius=1.5];
\graph [ empty nodes, nodes={xshift=5cm,circle,draw,fill=red,inner sep=2pt}, edges={thick}]  { subgraph I_n [n=5, clockwise,radius=15mm,V={a,b,c,d,e},phase=54]; {a,b,d} -> g[fill=white,shift=(18:.7)], {c,f[fill=white,shift=(108:1.55)]} -> d, g -> f -> e };
\foreach \x in {1,...,5} \node at ($ (5,0) + (54+72-72*\x:1.85) $) {$q_{\x}$};
\node at (6.05,.7) {\small $E_1$};
\node at (5.9,-.3) {\small $E_2$};
\node at (5.2,.55) {\small $E_4$};
\node at (4.9,-.3) {\small $E_4$};
\node at (4.5,-.8) {\small $E_3$};
\node at (4.65,.9) {\small $E_1$};
\node at (3.85,.2) {\small $E_5$};
\end{tikzpicture}
\caption{A weakly connected standard planar scaffold and its dual}
\label{example}
\end{figure}

In order to discuss duality of scaffolds whose edge weights are in the Bose--Mesner algebra $\bm{A}$, we consider the following situation.
Let $\mathsf{S}=\mathsf{S}(G,R;w)$ be a scaffold.
Suppose that $R$ is an ordered \emph{set}, i.e., the roots $r_1,\dots,r_{\ell}$ are mutually distinct, and that $G$ is weakly connected\footnote{Recall that $G$ is called \emph{weakly connected} if the underlying undirected graph of $G$ is connected. We may remark that every scaffold is a sum of scaffolds whose digraphs are weakly connected. Indeed, we can always add arbitrarily many edges with weight $J$ to a scaffold and then expand it using (AS\ref{AS2}).} and is given a fixed embedding (without crossings) in a closed disk in the plane, where the roots lie on the boundary circle, and the rest of the nodes as well as all the edges are in the interior of the disk.
We moreover assume that the roots are placed in clockwise order.\footnote{This assumption is more or less for ease of exposition when we consider a single scaffold and was not imposed in \cite{Martin2021LAA}, but it is essential in Conjecture \ref{Martin conjecture} below which involves multiple scaffolds. See Example \ref{standard placement}.}
In this case, let us call $\mathsf{S}$ a \emph{weakly connected standard planar scaffold}.
See the left diagram in Figure \ref{example}.

In the above situation, let $F$ be the set of faces of $G$ in the disk, i.e., the connected components of the region obtained by removing the edges and nodes of $G$ from the interior of the disk.
Suppose first that $\ell>0$.
Since $G$ is weakly connected, there exist exactly $\ell$ faces whose boundaries intersect with the boundary circle in arcs.
We call them $q_1,\dots,q_{\ell}$, where $q_i$ is incident with $r_i$ and $r_{i+1}$.
(We set $r_{\ell+1}:=r_1$.)
We define the dual digraph $G^{\dagger}$ of $G$ as the digraph having $F$ as the set of nodes, and whose edge set is $E^{\dagger}=\{e^{\dagger}:e\in E\}$, where $e^{\dagger}$ joins the two faces incident with $e$ and its direction is given by rotating $e$ through $90^{\circ}$ clockwise.\footnote{The direction of $e^{\dagger}$ is chosen counterclockwise in \cite{Martin2021LAA}, but we follow \cite[Section 10.2]{BM2008B} in this paper.}
We assume for simplicity that $w(e)\in \{A_0,\dots,A_d\}$ $(e\in E)$.
The \emph{dual scaffold} of $\mathsf{S}$ is the scaffold $\mathsf{S}^{\dagger}=(G^{\dagger},Q;w^{\dagger})$, where $Q=\{q_1,\dots,q_{\ell}\}$, and where $w^{\dagger}(e^{\dagger}):=E_i$ if $w(e)=A_i$ $(e\in E)$.
We plot the points representing $q_1,\dots,q_{\ell}$ on the incident arcs on the boundary circle, so that $\mathsf{S}^{\dagger}$ is again a weakly connected standard planar scaffold.\footnote{The dual scaffold $\mathsf{S}^{\dagger}$ appears to depend on the fixed embedding of $G$. That it is independent of the embedding is a special case of Conjecture \ref{Martin conjecture}, and is true for the class of translation association schemes by Theorem \ref{main theorem}.}
See Figure \ref{example}.
We proceed in the same manner when $\ell=0$, but we treat the unique face incident with the boundary circle as a non-root, and set $Q=\emptyset$.
In fact, the disk and its boundary circle play no essential role in this case; we may forget about them and view the unique unbounded face of $G$ as forming a non-root of $G^{\dagger}$.

In Conjecture \ref{Martin conjecture} below, which is a slight modification of \cite[Conjecture 4.1]{Martin2021LAA}, we will also view $\mathsf{S}$ and $\mathsf{S}^{\dagger}$ as functions of $(X,\mathcal{R})$.
In other words, we are given $G,R$, and also $\lambda\in \{0,1,\dots,c\}^E$, and then consider the scaffold $\mathsf{S}=\mathsf{S}(G,R;w)$ as well as its dual by setting $w(e):=A_{e^{\lambda}}$ $(e\in E)$ for the association scheme in question (which must have at least $c$ classes).

\begin{conj}[{cf.~\cite[Conjecture 4.1]{Martin2021LAA}}]\label{Martin conjecture}
Let $\{\mathsf{S}_{j,i}\}_{i=1}^{k_j}$ $(1\leqslant j\leqslant h)$ and $\{\mathsf{S}_i'\}_{i=1}^{k'}$ be sequences of weakly connected standard planar scaffolds whose edge weights are in $\{A_0,\dots,A_c\}$, such that the scaffolds in each sequence have the same order.
Assume that the following implication holds for all association schemes with $d\geqslant c$ classes:
\begin{equation*}
	\sum_{i=1}^{k_j} a_{j,i}\mathsf{S}_{j,i} =_j0 \quad (1\leqslant j\leqslant h) \quad \Longrightarrow \quad \sum_{i=1}^{k'} b_i \mathsf{S}_i'='0,
\end{equation*}
where $=_1,\dots,=_h,=' \,\in\{=,\ne\}$, and $a_{j,i},b_i\in\mathbb{C}$.
Then the following implication holds for all association schemes with $d\geqslant c$ classes:
\begin{equation*}
	\sum_{i=1}^{k_j} a_{j,i}|X|^{n_{j,i}}\mathsf{S}_{j,i}^{\dagger} =_j0 \quad (1\leqslant j\leqslant h) \quad \Longrightarrow \quad \sum_{i=1}^{k'} b_i |X|^{n_i'}\mathsf{S}_i^{\prime\dagger}='0,
\end{equation*}
where $n_{j,i}$ (resp.~$n_i'$) denotes the number of nodes of $\mathsf{S}_{j,i}$ (resp.~$\mathsf{S}_i'$).
\end{conj}

We have some comments on Conjecture \ref{Martin conjecture}.
First, while only equations were considered in \cite[Conjecture 4.1]{Martin2021LAA}, we allow inequations as well, because this makes it easier for example to interpret the $P$-polynomial and $Q$-polynomial properties.
We also allow that the scaffolds in different sequences have different orders.
Second, we dropped the condition in \cite[Conjecture 4.1]{Martin2021LAA} about bijections among the sets of root nodes of the scaffolds.
This condition seems unnecessary because the roots are already ordered by definition.
On the other hand, we added the condition on the clockwise placement of the root nodes on the boundary circle.
Third, and most importantly, we inserted the scalars $|X|^{n_{j,i}}$ and $|X|^{n_i'}$ in the dual (in)equations, which were missing in \cite[Conjecture 4.1]{Martin2021LAA}.
Theorem \ref{main theorem} below shows that Conjecture \ref{Martin conjecture} does not hold without this modification.
See also Examples \ref{necessity of scalars}, \ref{P- and/or Q-polynomial}, and \ref{standard placement}.

Suppose now that $X$ is endowed with the structure of an abelian group (written multiplicatively) with identity element $1$.
We call $(X,\mathcal{R})$ a \emph{translation} association scheme whenever $(xz,yz)\in R_i$ for all $0\leqslant i\leqslant d$, $(x,y)\in R_i$, and $z\in X$.
In this case, we can define the \emph{dual} $(X^*,\mathcal{R}^*)$ of $(X,\mathcal{R})$ on the character group $X^*$ of $X$; see Section \ref{sec: proof}.
The dual $(X^*,\mathcal{R}^*)$ is again a translation association scheme and has first and second eigenmatrices $P^*=Q$ and $Q^*=P$, respectively.
We now state our main result:

\begin{thm}\label{main theorem}
Conjecture \ref{Martin conjecture} holds true if we replace the phrases ``for all association schemes'' (two places) by ``for all translation association schemes''.
\end{thm}

Theorem \ref{main theorem} is immediate from

\begin{thm}\label{real main theorem}
Suppose that $(X,\mathcal{R})$ is a translation association scheme, and let $\mathsf{P}(\ell)$ (resp.~$\mathsf{P}^*(\ell)$) be the set of weakly connected standard planar scaffolds of order $\ell$ whose edge weights are in $\{A_0,\dots,A_d\}$ (resp.~$\{E_0^*,\dots,E_d^*\}$), where the $E_i^*$ are the primitive idempotents of $(X^*,\mathcal{R}^*)$.
When $\ell>0$, there is an isometric linear isomorphism from $\operatorname{span}_{\mathbb{C}}\mathsf{P}(\ell)$ to $\operatorname{span}_{\mathbb{C}}\mathsf{P}^*(\ell)$ which maps $\mathsf{S}\in\mathsf{P}(\ell)$ to $|X|^{n-\ell/2-1/2}\,\mathsf{S}^{\dagger}\in \mathsf{P}^*(\ell)$, where $\mathsf{S}^{\dagger}$ is taken\footnote{Here we are using the same convention as in the paragraph preceding Conjecture \ref{Martin conjecture}.} with respect to $(X^*,\mathcal{R}^*)$, and $n$ is the number of nodes in $\mathsf{S}$.
When $\ell=0$, with the same notation we instead have $\mathsf{S}=|X|^{n-1}\,\mathsf{S}^{\dagger}$ for $\mathsf{S}\in\mathsf{P}(0)$.
\end{thm}

\noindent
Theorem \ref{real main theorem} generalizes \cite[Proposition 11]{Jaeger1995JAC}, which deals with the case where $\ell=0$ and $d=|X|-1$.
We will prove Theorem \ref{real main theorem} in the next section.

We end this section with three additional examples which support our modification of \cite[Conjecture 4.1]{Martin2021LAA}.

\begin{exam}\label{necessity of scalars}
This example again illustrates the necessity of the scalars $|X|^{n_{j,i}}$ and $|X|^{n_i'}$ in the dual (in)equations in Conjecture \ref{Martin conjecture}.
We have
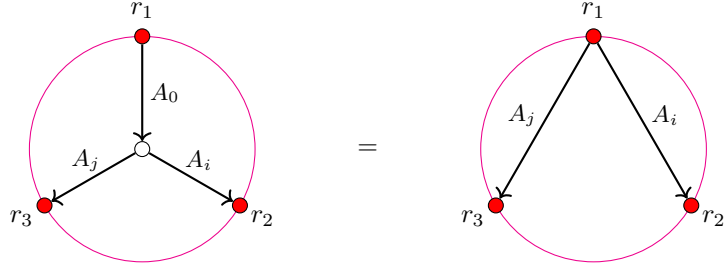
\begin{center}
\begin{tikzpicture}
\draw[magenta] (0,0) circle [radius=1.5];
\graph [ empty nodes, nodes={circle,draw,fill=red,inner sep=2pt}, edges={thick}]  { subgraph I_n [n=3, clockwise,radius=15mm,V={a,b,c}];  a -> d[fill=white] -> {b,c} };
\foreach \x in {1,2,3} \node at (90+120-120*\x:1.85) {$r_{\x}$};
\node at (.3,.75) {\small $A_0$};
\node at (.75,-.15) {\small $A_i$};
\node at (-.75,-.15) {\small $A_j$};
\node at (3,0) {$=$};
\draw[magenta] (6,0) circle [radius=1.5];
\graph [ empty nodes, nodes={xshift=6cm,circle,draw,fill=red,inner sep=2pt}, edges={thick}]  { subgraph I_n [n=3, clockwise,radius=15mm,V={a,b,c}]; a -> {b,c} };
\foreach \x in {1,2,3} \node at ($ (6,0) + (90+120-120*\x:1.85) $) {$r_{\x}$};
\node at (6.95,.45) {\small $A_i$};
\node at (5.05,.45) {\small $A_j$};
\end{tikzpicture}
\end{center}
since $A_0=I$.
The LHS (resp.~RHS) above has four (resp.~three) nodes.
Hence the dual equation, divided by $|X|^3$, becomes
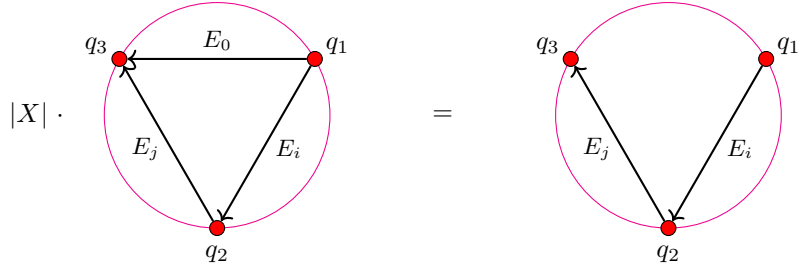
\begin{center}
\begin{tikzpicture}
\node at (-2.4,0) {$|X|\ \cdot$};
\draw[magenta] (0,0) circle [radius=1.5];
\graph [ empty nodes, nodes={circle,draw,fill=red,inner sep=2pt}, edges={thick}]  { subgraph I_n [n=3, clockwise,radius=15mm,V={a,b,c},phase=30];  a -> {b,c}, b -> c };
\foreach \x in {1,2,3} \node at (30+120-120*\x:1.85) {$q_{\x}$};
\node at (0,1) {\small $E_0$};
\node at (.95,-.45) {\small $E_i$};
\node at (-.95,-.45) {\small $E_j$};
\node at (3,0) {$=$};
\draw[magenta] (6,0) circle [radius=1.5];
\graph [ empty nodes, nodes={xshift=6cm,circle,draw,fill=red,inner sep=2pt}, edges={thick}]  { subgraph I_n [n=3, clockwise,radius=15mm,V={a,b,c},phase=30]; a -> b -> c };
\foreach \x in {1,2,3} \node at ($ (6,0) + (30+120-120*\x:1.85) $) {$q_{\x}$};
\node at (6.95,-.45) {\small $E_i$};
\node at (5.05,-.45) {\small $E_j$};
\end{tikzpicture}
\end{center}
which is also immediate from $|X|E_0=J$.
These are special cases of Rules \textsf{SR}$0$ and \textsf{SR}$0'$ in \cite[Appendix A]{Martin2021LAA}.
\end{exam}

\begin{exam}\label{P- and/or Q-polynomial}
Recall that the \emph{intersection numbers} $p_{i,j}^k$ and the \emph{Krein parameters} $q_{i,j}^k$ are defined respectively by
\begin{equation*}
	A_iA_j=\sum_{k=0}^d p_{i,j}^k A_k, \qquad E_i\circ E_j=|X|^{-1}\sum_{k=0}^d q_{i,j}^k E_k \qquad (0\leqslant i,j\leqslant d),
\end{equation*}
and that $p_{i,j}^k=0$ (resp.~$q_{i,j}^k=0$) if and only if the scaffold in \eqref{triangle scaffold} (resp.~\eqref{star scaffold}) equals zero.
See, e.g., \cite[p.~128]{BI1984B} or \cite[Lemma 3.2]{Terwilliger1992JAC} for the latter.
Suppose now that $(X,\mathcal{R})$ is $P$-\emph{polynomial} (or \emph{metric}), i.e., that the $A_i$ are symmetric matrices and that, for all $0\leqslant i,j,k\leqslant d$, we have $p_{i,j}^k=0$ (resp.~$p_{i,j}^k\ne 0$) whenever one of $i,j,k$ is greater than (resp.~equal to) the sum of the other two.
Note that the former condition translates to
\begin{center}
\begin{tikzpicture}
\graph [ empty nodes, nodes={circle,draw,fill=red,inner sep=2pt}, edges={thick}] { a[label=left:$r_1$] ->["\small $A_i$"] b[label=right:$r_2$] };
\node at (2.5,0) {$=$};
\graph [ empty nodes, nodes={xshift=3.9cm,circle,draw,fill=red,inner sep=2pt}, edges={thick}]
{ a[label=left:$r_1$] <-["\small $A_i$"] b[label=right:$r_2$] };
\node at (7.2,0) {$(0\leqslant i\leqslant d)$,};
\end{tikzpicture}
\end{center}
and that the latter is a set of equations and inequations on scaffolds in \eqref{triangle scaffold}.
In this case, we for example have
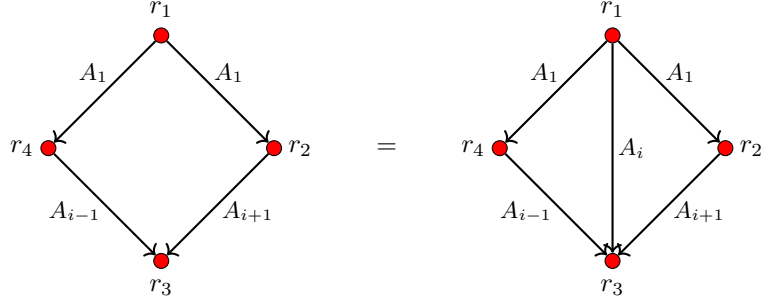
\begin{center}
\begin{tikzpicture}
\graph [ empty nodes, nodes={circle,draw,fill=red,inner sep=2pt}, edges={thick}]  { subgraph I_n [n=4, clockwise,radius=15mm,V={a,b,c,d}];  a -> {b,d} -> c };
\foreach \x in {1,2,3,4} \node at (90+90-90*\x:1.85) {$r_{\x}$};
\node at (.9,1) {\small $A_1$};
\node at (1.15,-.85) {\small $A_{i+1}$};
\node at (-.9,1) {\small $A_1$};
\node at (-1.15,-.85) {\small $A_{i-1}$};
\node at (3,0) {$=$};
\graph [ empty nodes, nodes={xshift=6cm,circle,draw,fill=red,inner sep=2pt}, edges={thick}]  { subgraph I_n [n=4, clockwise,radius=15mm,V={a,b,c,d}]; a -> {b,d} -> c, a -> c };
\foreach \x in {1,2,3,4} \node at ($ (6,0) + (90+90-90*\x:1.85) $) {$r_{\x}$};
\node at (6.9,1) {\small $A_1$};
\node at (7.15,-.85) {\small $A_{i+1}$};
\node at (5.1,1) {\small $A_1$};
\node at (4.85,-.85) {\small $A_{i-1}$};
\node at (6.25,0) {\small $A_i$};
\end{tikzpicture}
\end{center}
for $1\leqslant i\leqslant d-1$.
This is a special case of Rule \textsf{SR}$4'$ in \cite[Appendix A]{Martin2021LAA}.
Note that we have used equations and inequations of scaffolds of order two, three, and four.
The $Q$-\emph{polynomial} (or \emph{cometric}) property is defined analogously by replacing the $A_i$ by the $E_i$ and the $p_{i,j}^k$ by the $q_{i,j}^k$.
The statement dual to the above follows from Rule \textsf{SR}$4$ in \cite[Appendix A]{Martin2021LAA}.
\end{exam}

\begin{exam}\label{standard placement}
This example shows that Conjecture \ref{Martin conjecture} fails without the condition on the clockwise placement of the roots.
It is clear that
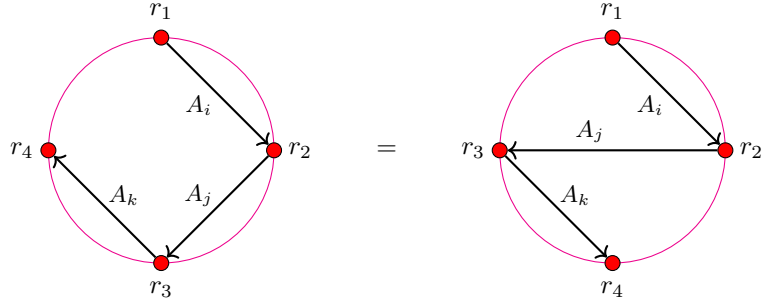
\begin{center}
\begin{tikzpicture}
\draw[magenta] (0,0) circle [radius=1.5];
\graph [ empty nodes, nodes={circle,draw,fill=red,inner sep=2pt}, edges={thick}]  { subgraph I_n [n=4, clockwise,radius=15mm,V={a,b,c,d}];  a -> b -> c -> d };
\foreach \x in {1,2,3,4} \node at (90+90-90*\x:1.85) {$r_{\x}$};
\node at (.5,.6) {\small $A_i$};
\node at (.5,-.6) {\small $A_j$};
\node at (-.5,-.6) {\small $A_k$};
\node at (3,0) {$=$};
\draw[magenta] (6,0) circle [radius=1.5];
\graph [ empty nodes, nodes={xshift=6cm,circle,draw,fill=red,inner sep=2pt}, edges={thick}]  { subgraph I_n [n=4, clockwise,radius=15mm,V={a,b,c,d}]; a -> b -> d -> c };
\node at ($(6,0)+(90:1.85)$) {$r_1$};
\node at ($(6,0)+(0:1.85)$) {$r_2$};
\node at ($(6,0)+(270:1.85)$) {$r_4$};
\node at ($(6,0)+(180:1.85)$) {$r_3$};
\node at (6.5,.6) {\small $A_i$};
\node at (5.7,.25) {\small $A_j$};
\node at (5.5,-.6) {\small $A_k$};
\end{tikzpicture}
\end{center}
where the LHS is a weakly connected standard planar scaffold, whereas the RHS is not.
If we take the dual of both sides (but without specifying a particular ordering of the faces for the RHS), we obtain
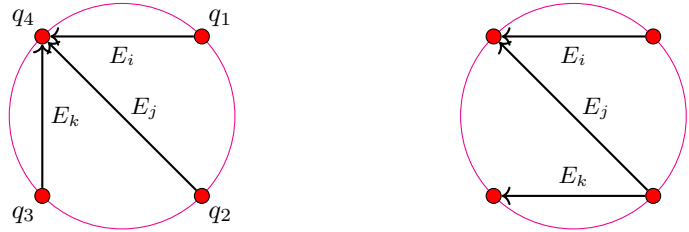
\begin{center}
\begin{tikzpicture}
\draw[magenta] (0,0) circle [radius=1.5];
\graph [ empty nodes, nodes={circle,draw,fill=red,inner sep=2pt}, edges={thick}]  { subgraph I_n [n=4, clockwise,phase=45,radius=15mm,V={a,b,c,d}];  { a,b,c } -> d };
\foreach \x in {1,2,3,4} \node at (45+90-90*\x:1.85) {$q_{\x}$};
\node at (0,.8) {\small $E_i$};
\node at (.3,.1) {\small $E_j$};
\node at (-.75,0) {\small $E_k$};
\draw[magenta] (6,0) circle [radius=1.5];
\graph [ empty nodes, nodes={xshift=6cm,circle,draw,fill=red,inner sep=2pt}, edges={thick}]  { subgraph I_n [n=4, clockwise,phase=45,radius=15mm,V={a,b,c,d}]; a -> d <- b -> c };
\node at (6,.8) {\small $E_i$};
\node at (6.3,.1) {\small $E_j$};
\node at (6,-.8) {\small $E_k$};
\end{tikzpicture}
\end{center}
respectively, but these are distinct in general, for example when $j=0$ and $i,k>0$.
\end{exam}

\section{Proof of Theorem \ref{real main theorem}}
\label{sec: proof}

We begin by recalling necessary facts from \cite[Section 6]{MT2009EJC} regarding the translation association scheme $(X,\mathcal{R})$.
Let $X^*$ be the character group of $X$ with identity $\iota$, i.e., $\iota(x)=1$ $(x\in X)$.
Let
\begin{equation*}
	\hat{\varepsilon}:=|X|^{-1/2}\sum_{x\in X}\overline{\varepsilon(x)}\hat{x} \qquad (\varepsilon\in X^*),
\end{equation*}
where $\bar{\ }$ denotes complex conjugate.
Then the $\hat{\varepsilon}$ form an orthonormal basis of $\mathbb{C}^X$, and we may identify $\mathbb{C}^{X^*}$ with $\mathbb{C}^X$ (and hence $\mathbb{C}^{(X^*)^{\ell}}$ with $\mathbb{C}^{X^{\ell}}$) via this basis.
Note that the (di)graphs $(X,R_i)$ $(1\leqslant i\leqslant d)$ are Cayley graphs, so that the $\hat{\varepsilon}$ are common eigenvectors of the Bose--Mesner algebra $\bm{A}$.
Hence there is a partition $X^*=X_0^*\sqcup\cdots\sqcup X_d^*$ such that $\{\hat{\varepsilon}:\varepsilon\in X_i^*\}$ is an orthonormal basis of $E_i\mathbb{C}^X$ $(0\leqslant i\leqslant d)$.
Let $\mathcal{R}^*=\{R_0^*,\dots,R_d^*\}$ be the set of binary relations on $X^*$ given by
\begin{equation}\label{dual relations}
	R_i^*=\{(\varepsilon,\eta)\in(X^*)^2:\eta\varepsilon^{-1}\in X_i^*\} \qquad (0\leqslant i\leqslant d).
\end{equation}
Then it follows that $(X^*,\mathcal{R}^*)$ is again a translation association scheme, called the \emph{dual} of $(X,\mathcal{R})$.
We mentioned earlier that it has first and second eigenmatrices $P^*=Q$ and $Q^*=P$, respectively.

We first prove Theorem \ref{real main theorem} when $\ell>0$.
For $\bm{\varepsilon}=(\varepsilon_1,\dots,\varepsilon_{\ell})\in (X^*)^{\ell}$, we write $\hat{\bm{\varepsilon}}=\hat{\varepsilon}_1\otimes\cdots\otimes\hat{\varepsilon}_{\ell}$ for brevity.
Let $\Gamma$ be the subgroup of $(X^*)^{\ell}$ defined by
\begin{equation*}
	\Gamma=\{(\varepsilon_1,\dots,\varepsilon_{\ell})\in(X^*)^{\ell}:\varepsilon_1\cdots\varepsilon_{\ell}=\iota\},
\end{equation*}
and let
\begin{equation*}
	\mathbb{C}^{\Gamma}=\operatorname{span}_{\mathbb{C}}\{\hat{\bm{\varepsilon}}:\bm{\varepsilon}\in\Gamma\} \subset \mathbb{C}^{X^{\ell}}.
\end{equation*}
Define the surjective homomorphism $\Phi:(X^*)^{\ell}\rightarrow \Gamma$ by
\begin{equation}\label{f}
	\Phi(\varepsilon_1,\dots,\varepsilon_{\ell})=(\varepsilon_{\ell}\varepsilon_1^{-1},\varepsilon_1\varepsilon_2^{-1},\dots,\varepsilon_{\ell-1}\varepsilon_{\ell}^{-1}) \qquad (\varepsilon_1,\dots,\varepsilon_{\ell}\in X^*),
\end{equation}
and let $\Psi:\mathbb{C}^{\Gamma}\rightarrow\mathbb{C}^{X^{\ell}}$ be the linear map defined by
\begin{equation*}
	\Psi(\hat{\bm{\varepsilon}})=|X|^{-1/2} \!\! \sum_{\substack{\bm{\eta}\in (X^*)^{\ell} \\ \Phi(\bm{\eta})=\bm{\varepsilon}}} \!\! \hat{\bm{\eta}} \qquad (\bm{\varepsilon}\in\Gamma).
\end{equation*}
Note that the sum on the RHS has $|X|$ terms, so that $\Psi$ is an isometry.

Let $\mathsf{S}=\mathsf{S}(G,R;w)\in\mathsf{P}(\ell)$, and let $\lambda\in\{0,\dots,d\}^E$ be such that $w(e)=A_{e^{\lambda}}$ $(e\in E)$.
Hence
\begin{equation*}
	\mathsf{S} = \sum_{\sigma\in X^V}\!\left(\prod_{e\in E}(A_{e^{\lambda}})_{t(e)^{\sigma}\!,h(e)^{\sigma}} \!\right) \widehat{\bm{r}^{\sigma}}.
\end{equation*}
We assume that $G$ has $n$ nodes and $m$ edges.

For $0\leqslant i\leqslant d$ and $x,y\in X$, by \eqref{eigenmatrices} we have
\begin{align*}
	(A_i)_{x,y} &= \langle A_i\hat{y},\hat{x}\rangle \\
	&= \sum_{\varepsilon\in X^*} \langle A_i\hat{\varepsilon}\overline{\hat{\varepsilon}}{}^{\,\mathsf{T}}\hat{y},\hat{x}\rangle \\
	&=\sum_{\varepsilon\in X^*} P_{\varepsilon, i} \langle \hat{\varepsilon}\overline{\hat{\varepsilon}}{}^{\,\mathsf{T}}\hat{y},\hat{x}\rangle \\
	&=\sum_{\varepsilon\in X^*} P_{\varepsilon, i} \langle \hat{\varepsilon},\hat{x}\rangle\langle\hat{y},\hat{\varepsilon}\rangle \\
	&=|X|^{-1}\sum_{\varepsilon\in X^*} P_{\varepsilon, i}\, \overline{\varepsilon(x)}\varepsilon(y),
\end{align*}
where we write
\begin{equation}\label{P}
	P_{\varepsilon, i}:=P_{j, i} \quad \text{if} \  \varepsilon\in X_j^*.
\end{equation}
Hence, by setting for $\tau\in (X^*)^E$,
\begin{equation}\label{boundary homomorphism}
	v^{\tau}= \Bigg(\prod_{\substack{e\in E \\ h(e)=v}} \!\! e^{\tau}\!\Bigg)^{\!\!\!-1} \! \Bigg(\prod_{\substack{e\in E \\ t(e)=v}} \!\! e^{\tau}\!\Bigg) \qquad (v\in V)
\end{equation}
and $\bm{r}^{\tau}=(r_1^{\tau},\dots,r_{\ell}^{\tau})$ for brevity, we have
\begin{align*}
	\mathsf{S} &= |X|^{-m} \! \sum_{\sigma\in X^V}\sum_{\tau\in (X^*)^E} \!\left( \prod_{e\in E}P_{e^{\tau}\!,e^{\lambda}} \overline{e^{\tau}(t(e)^{\sigma})}e^{\tau}(h(e)^{\sigma})\!\right) \widehat{\bm{r}^{\sigma}} \\
	& = |X|^{-m} \!\! \sum_{\tau\in(X^*)^E} \!\left( \prod_{e\in E}P_{e^{\tau}\!,e^{\lambda}}\!\right) \! \sum_{\sigma\in X^V} \!\Bigg(\prod_{v\in V} \overline{v^{\tau}(v^{\sigma})}\Bigg)\, \widehat{\bm{r}^{\sigma}} \\
	&= |X|^{-m} \!\! \sum_{\tau\in(X^*)^E} \!\left( \prod_{e\in E}P_{e^{\tau}\!,e^{\lambda}}\!\right) \!\!\Bigg(\prod_{v\in V\setminus R} \!\! |X|\delta_{v^{\tau}\!,\iota}\!\Bigg)\, |X|^{\ell/2} \,\widehat{\bm{r}^{\tau}} \\
	&= |X|^{n-m-\ell/2} \!\! \sum_{\tau\in (X^*)^E} \!\left( \prod_{e\in E}P_{e^{\tau}\!,e^{\lambda}}\!\right)\!\!\Bigg(\prod_{v\in V\setminus R} \!\! \delta_{v^{\tau}\!,\iota}\!\Bigg)\,\widehat{\bm{r}^{\tau}}.
\end{align*}
In particular, it follows that $\mathsf{S}\in \mathbb{C}^{\Gamma}$ since
\begin{equation*}
	r_1^{\tau}\cdots r_{\ell}^{\tau}\prod_{v\in V\setminus R} \!\! v^{\tau} = \prod_{e\in E} e^{\tau}(e^{\tau})^{-1} = \iota,
\end{equation*}
and hence that $\operatorname{span}_{\mathbb{C}}\mathsf{P}(\ell)\subset\mathbb{C}^{\Gamma}$.

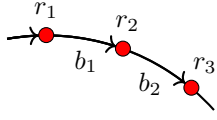
\begin{figure}
\begin{tikzpicture}[node/.style={circle,fill=red,draw=black,inner sep=2pt}]
\draw [->, thick] (0,0) arc [start angle=100, end angle=92, radius=3cm];
\draw [->, thick] (0,0) arc [start angle=100, end angle=72, radius=3cm];
\draw [->, thick] (0,0) arc [start angle=100, end angle=52, radius=3cm];
\draw [thick] (0,0) arc [start angle=100, end angle=42, radius=3cm];
\node [node,label=above:$r_1$] at (.52,.05) {};
\node [node,label=above:$\,\,r_2$] at (1.54,-.13) {};
\node [node,label=above right:$\!\!\!r_3$] at (2.45,-.65) {};
\node at (1.05,-.3) {$b_1$};
\node at (1.9,-.6) {$b_2$};
\end{tikzpicture}
\caption{Boundary arcs}
\label{boundary arcs}
\end{figure}

Now we have
\begin{equation*}
	\Psi(\mathsf{S}) =  |X|^{n-m-\ell/2-1/2} \!\! \sum_{\tau\in (X^*)^E} \!\left( \prod_{e\in E}P_{e^{\tau}\!,e^{\lambda}}\!\right)\!\!\Bigg(\prod_{v\in V\setminus R} \!\! \delta_{v^{\tau}\!,\iota}\!\Bigg)\!\sum_{\substack{\bm{\eta}\in (X^*)^{\ell} \\ \Phi(\bm{\eta})=\bm{r}^{\tau} }} \!\! \hat{\bm{\eta}}.
\end{equation*}
Let $B=\{b_1,\dots,b_{\ell}\}$ be the set of arcs obtained by removing the root nodes from the boundary circle, where $b_i$ joins $r_i$ and $r_{i+1}$.
(Recall that $r_{\ell+1}:=r_1$.)
We give an orientation on these arcs clockwise, i.e., $t(b_i)=r_i$ and $h(b_i)=r_{i+1}$.
See Figure \ref{boundary arcs}.
To every pair $(\tau,\bm{\eta})$ with $\tau\in (X^*)^E$ and $\bm{\eta}=(\eta_1,\dots,\eta_{\ell})\in (X^*)^{\ell}$, we associate $\tilde{\tau}\in (X^*)^{E\sqcup B}$ as follows:
\begin{equation*}
	e^{\tilde{\tau}}=\begin{cases} e^{\tau} & \text{if} \ e\in E, \\ \eta_i & \text{if} \ e=b_i \ \text{for some} \ i. \end{cases}
\end{equation*}
We also define $v^{\tilde{\tau}}$ $(v\in V)$ in the same manner as \eqref{boundary homomorphism}, taking into account the arcs in $B$ as well.
Then, in view of the definition of $\Phi$ (cf.~\eqref{f}), that $\Phi(\bm{\eta})=\bm{r}^{\tau}$ is equivalent to saying that $(r_i)^{\tilde{\tau}}=\iota$ for $1\leqslant i\leqslant\ell$.
Hence we may rewrite $\Psi(\mathsf{S})$ as
\begin{equation*}
	\Psi(\mathsf{S}) =  |X|^{n-m-\ell/2-1/2} \sum_{(\tau,\bm{\eta})} \!\left( \prod_{e\in E}P_{e^{\tau}\!,e^{\lambda}}\!\right)\! \hat{\bm{\eta}},
\end{equation*}
where the sum is over $(\tau,\bm{\eta})\in (X^*)^E\times (X^*)^{\ell}$ such that $v^{\tilde{\tau}}=\iota$ for all $v\in V$; that is to say, $\tilde{\tau}$ is an $X^*$-\emph{circulation} on the plane digraph $\tilde{G}:=(V,E\sqcup B)$.

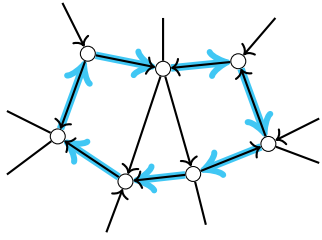
\begin{figure}
\begin{tikzpicture}[new set=import nodes]
\begin{scope}[nodes={set=import nodes}, node/.style={circle,fill=white,draw=black,inner sep=2pt}]
\node[node] (a) at (0,0) {};
\node[node] (b) at (1,-.2) {};
\node[node] (c) at (2,-.1) {};
\node[node] (d) at (2.4,-1.2) {};
\node[node] (e) at (1.4,-1.6) {};
\node[node] (f) at (.5,-1.7) {};
\node[node] (g) at (-.4,-1.1) {};
\node (a1) at ($ (a) + (-.4,.8) $) {};
\node (b1) at ($ (b) + (0,.8) $) {};
\node (c1) at ($ (c) + (.6,.7) $) {};
\node (d1) at ($ (d) + (.8,.4) $) {};
\node (d2) at ($ (d) + (.8,-.3) $) {};
\node (e1) at ($ (e) + (.2,-.8) $) {};
\node (f1) at ($ (f) + (-.3,-.8) $) {};
\node (g1) at ($ (g) + (-.8,.4) $) {};
\node (g2) at ($ (g) + (-.8,-.6) $) {};
\end{scope}
\graph[edges={thick,line width=1mm,draw=cyan!60,>={Classical TikZ Rightarrow[width=4mm,length=3mm]},shorten >=-0pt}] {
(import nodes);
a -> b -> c -> d -> e -> f -> g -> a };
\graph[edges={thick}] {
(import nodes);
a -> b <- c <- d <- e -> f -> g <- a,
b -> { e,f },
a1 -> a,
b1 -- b,
c1 -> c,
d1 -> d,
d2 -- d,
e1 -- e,
f1 -> f,
{ g1,g2 } -- g,
};
\end{tikzpicture}
\caption{A cycle-orientation of a cycle}
\label{cycle}
\end{figure}

Examples of $X^*$-circulations are obtained from cycles in the underlying graph of $\tilde{G}$.
Suppose that $C\subset E\sqcup B$ forms a cycle (when the directions are discarded), and pick one of the two cycle-orientations of $C$.
See Figure \ref{cycle}.
Let $\varepsilon\in X^*$, and define $\alpha=\alpha_{C,\varepsilon}\in (X^*)^{E\sqcup B}$ as follows:
for $e\in E\sqcup B$, we set $e^{\alpha}=\varepsilon$ if $e\in C$ and its direction agrees with the chosen cycle-orientation, $e^{\alpha}=\varepsilon^{-1}$ if $e\in C$ and its direction is the opposite of the cycle-orientation, and $e^{\alpha}=\iota$ if $e\not\in C$.
Then $\alpha$ is an $X^*$-circulation on $\tilde{G}$.
Moreover, it is well known\footnote{See, e.g., \cite[Proposition 7.13]{BM2008B} for the case of $\mathbb{R}$-circulations. If $\alpha$ is a nontrivial $X^*$-circulation, then its support must contain a cycle $C$, and we can choose $\varepsilon\in X^*$ so that $\beta=\alpha\cdot\alpha_{C,\varepsilon}^{-1}$ has strictly smaller support than $\alpha$.} that the subgroup $Z=Z(\tilde{G},X^*)$ of $(X^*)^{E\sqcup B}$ consisting of the $X^*$-circulations on $\tilde{G}$ is generated by the elements of this form.
Since $\tilde{G}$ is a plane digraph, it suffices to take only the $\alpha_{C,\varepsilon}$ where $C$ is the set of boundary edges of a face of $\tilde{G}$.

Let $\infty$ denote the unbounded face of $\tilde{G}$, i.e., the complement of the disk in the plane, so that the set of faces of $\tilde{G}$ is $F\sqcup\{\infty\}$.
For every $f\in F\sqcup\{\infty\}$, let $C_f$ be the set of boundary edges of $f$, and choose the cycle-orientation so that $f$ is on the right.
Then we have a surjective homomorphism $\Theta:(X^*)^{F\sqcup\{\infty\}}\rightarrow Z$ given by
\begin{equation*}
	\Theta(\mu) =\prod_{f\in F\sqcup\{\infty\}} \!\!\!\! \alpha_{C_f,f^{\mu}} \qquad (\mu\in (X^*)^{F\sqcup\{\infty\}}).
\end{equation*}
More specifically,
\begin{equation*}
	e^{\Theta(\mu)}=h(e^{\dagger})^{\mu} (t(e^{\dagger})^{\mu})^{-1} \qquad (e\in E\sqcup B).
\end{equation*}
Note that $\Theta$ is a $|X^*|$-to-one function, and it becomes an isomorphism if we restrict it to the subgroup consisting of those $\mu\in (X^*)^{F\sqcup\{\infty\}}$ such that $\infty^{\mu}=\iota$.
Note also that for such $\mu$ we have
\begin{equation*}
	b_i^{\Theta(\mu)}=q_i^{\mu} \qquad (1\leqslant i\leqslant \ell),
\end{equation*}
since $h(b_i^{\dagger})=q_i$ and $t(b_i^{\dagger})=\infty$.

Combining these comments, we have
\begin{align*}
	\Psi(\mathsf{S}) &= |X|^{n-m-\ell/2-1/2} \!\!\!\! \sum_{\substack{\mu\in (X^*)^{F\sqcup\{\infty\}} \\ \infty^{\mu}=\iota }} \!\left( \prod_{e\in E}P_{e^{\Theta(\mu)}\!,e^{\lambda}}\!\right) \widehat{\bm{b}^{\Theta(\mu)}} \\
	&= |X|^{n-m-\ell/2-1/2} \!\! \sum_{\mu\in (X^*)^F} \!\left( \prod_{e\in E}P_{h(e^{\dagger})^{\mu} (t(e^{\dagger})^{\mu})^{-1}\!,e^{\lambda}}\!\right) \widehat{\bm{q}^{\mu}},
\end{align*}
where $\bm{b}^{\Theta(\mu)}=(b_1^{\Theta(\mu)},\dots,b_{\ell}^{\Theta(\mu)})$ and $\bm{q}^{\mu}=(q_1^{\mu},\dots,q_{\ell}^{\mu})$.
In the last expression above, assume that $h(e^{\dagger})^{\mu} (t(e^{\dagger})^{\mu})^{-1}\in X_j^*$, or equivalently, $(t(e^{\dagger})^{\mu},h(e^{\dagger})^{\mu})\in R_j^*$ (cf.~\eqref{dual relations}).
Then, by \eqref{eigenmatrices} (applied to $(X^*,\mathcal{R}^*)$), \eqref{P}, and $P=Q^*$, we have
\begin{equation*}
	P_{h(e^{\dagger})^{\mu} (t(e^{\dagger})^{\mu})^{-1}\!,e^{\lambda}} = P_{j,e^{\lambda}} = Q_{j,e^{\lambda}}^* = |X^*| (E_{e^{\lambda}}^*)_{t(e^{\dagger})^{\mu}\!,h(e^{\dagger})^{\mu}}.
\end{equation*}
It follows that
\begin{align*}
	\Psi(\mathsf{S}) &= |X|^{n-\ell/2-1/2} \!\! \sum_{\mu\in (X^*)^F} \!\left( \prod_{e\in E}(E_{e^{\lambda}}^*)_{t(e^{\dagger})^{\mu}\!,h(e^{\dagger})^{\mu}}\!\right) \widehat{\bm{q}^{\mu}} \\
	&= |X|^{n-\ell/2-1/2}\,\mathsf{S}^{\dagger},
\end{align*}
where $\mathsf{S}^{\dagger}$ is taken with respect to the dual $(X^*,\mathcal{R}^*)$.

We have shown that $\operatorname{span}_{\mathbb{C}}\mathsf{P}(\ell)\subset \mathbb{C}^{\Gamma}$, and that $\Psi$ gives rise to a linear isometry from $\operatorname{span}_{\mathbb{C}}\mathsf{P}(\ell)$ to $\operatorname{span}_{\mathbb{C}}\mathsf{P}^*(\ell)$.
That it is surjective and hence is an isomorphism, is immediate from $G^{\dagger\dagger\dagger\dagger}=G$, where we recall that $G$ is weakly connected.
This establishes Theorem \ref{real main theorem} when $\ell>0$.

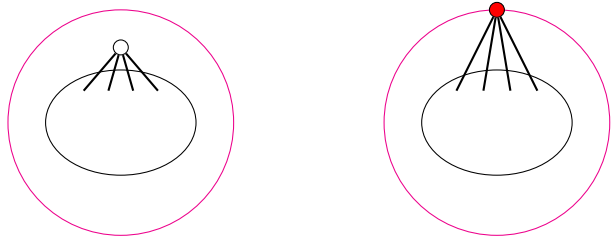
\begin{figure}
\begin{tikzpicture}[new set=import nodes]
\draw[magenta] (0,0) circle [radius=1.5];
\draw[magenta] (5,0) circle [radius=1.5];
\begin{scope}[nodes={set=import nodes}, node/.style={circle,draw=black,inner sep=2pt}]
\node[node,fill=white] (a) at (0,1) {};
\node (b) at (-.6,.3) {};
\node (c) at (-.2,.3) {};
\node (d) at (.2,.3) {};
\node (e) at (.6,.3) {};
\node[node,fill=red] (f) at ($ (a) + (5,.5) $) {};
\node (g) at ($ (b) + (5,0) $) {};
\node (h) at ($ (c) + (5,0) $) {};
\node (i) at ($ (d) + (5,0) $) {};
\node (j) at ($ (e) + (5,0) $) {};
\end{scope}
\draw (0,0) circle [x radius=1, y radius=.7];
\draw (5,0) circle [x radius=1, y radius=.7];
\graph[edges={thick}] {
(import nodes);
a -- { b,c,d,e }, f -- { g,h,i,j }
};
\end{tikzpicture}
\caption{An $\mathsf{S}$ with $\ell=0$ and its deformation $\tilde{\mathsf{S}}$}
\label{deformation}
\end{figure}

It remains to handle the case where $\ell=0$.
It is possible and is not difficult to adjust the above discussions to this case, but we may also argue as follows.
Pick any node which is incident with the outmost face of $G$ in the disk, shift it to the boundary circle, and change it to a root.
Call the resulting scaffold of order one $\tilde{\mathsf{S}}$.
See Figure \ref{deformation}.
Let $\mathsf{S}_1$ be the scaffold with diagram %
\begin{tikzpicture}
\graph [ empty nodes, nodes={circle,draw,fill=red,inner sep=2pt}] { a };
\end{tikzpicture},
and observe that $\mathsf{S}=\langle \tilde{\mathsf{S}}, \mathsf{S}_1 \rangle$.
Observe also that $\tilde{\mathsf{S}}^{\dagger}$ is obtained from $\mathsf{S}^{\dagger}$ by applying the same procedure to the outmost face of $G$ (as a node of $G^{\dagger}$).
Since $\mathsf{S}_1^{\dagger}$ again has diagram 
\begin{tikzpicture}
\graph [ empty nodes, nodes={circle,draw,fill=red,inner sep=2pt}] { a };
\end{tikzpicture}
(but with respect to $(X^*,\mathcal{R}^*)$), we have $\mathsf{S}^{\dagger}=\langle \tilde{\mathsf{S}}^{\dagger}, \mathsf{S}_1^{\dagger} \rangle$.
Hence it follows that
\begin{equation*}
	\mathsf{S}=\langle \tilde{\mathsf{S}}, \mathsf{S}_1 \rangle = \langle \Psi(\tilde{\mathsf{S}}), \Psi(\mathsf{S}_1) \rangle =\langle |X|^{n-1}\tilde{\mathsf{S}}^{\dagger}, \mathsf{S}_1^{\dagger} \rangle = |X|^{n-1} \mathsf{S}^{\dagger},
\end{equation*}
as desired.

This completes the proof of Theorem \ref{real main theorem}.

\section*{Acknowledgments}

The authors thank the anonymous referee for offering many valuable suggestions and comments.




\end{document}